\documentclass[12pt, a4paper]{article}
\usepackage{}
\usepackage{bbm}
\usepackage{mathrsfs}
\usepackage{amsfonts}
\usepackage{mathrsfs}

\usepackage{latexsym}
\usepackage{graphicx}
\usepackage{amssymb}

\newtheorem{theorem}{Theorem}[section]
\newtheorem{lemma}[theorem]{Lemma}
\newtheorem{corollary}[theorem]{Corollary}

\newtheorem{definition}[theorem]{Definition}

\def\gcd{{\rm gcd}}

\title{{\bf Exponents of the primitive Boolean matrices with fixed girth}\thanks{supported by Jiangsu students' training program for innovation and entrepreneurship (No. 201410324005Z, national level and provincial level), Jiangsu Qing Lan Project (2014A), JSNSFC (BK2012245), NSFC
(Nos. 11271315, 11171290, 11201417).}}
\author{Guanglong Yu\thanks{E-mail addresses:
yglong01@163.com.} ~ Shuguang Guo ~ Wenjuan Jia ~ Yuhan Che~
\\ ~  \\ {\footnotesize
Department of Mathematics and statistics, Yancheng Teachers University,}\\
{\footnotesize  Yancheng, 224002, Jiangsu, P.R. China}}

\date{}

\begin{document}
\maketitle

\begin{abstract}
The $girth$ of a primitive Boolean matrix is defined to be the $girth$ of its associated digraph. In this paper, among all primitive Boolean matrices of order $n$, the primitive exponents of those of girth $g$ are considered. For the primitive matrices of both order $n\geq 10$ and girth $g>\frac{n^{2}-4n}{4(n-3)}$, the matrices with primitive exponents in $[2n-2 +(g- 1)(n-3), n+g(n-2)]$ are completely characterized.


\end{abstract}

\section{Introduction}

\ \ \ \ We adopt the standard conventions, notations and definitions for
Boolean matrices, their entries,
arithmetics, powers and primitivity. The reader who is not familiar with these
matters is referred to \cite{BLiu}--\cite{LB1} and \cite{JYS}.

In this paper, we permit loops but no multiple arcs in a
digraph. For a digraph $S$, we denote by $V(S)$ the vertex set and denoted by $E(S)$ the
arc set.

\begin{definition} \label{de1.2.2}
Let $A$ be a square $(0,1)$-Boolean matrix of order $n$. The {\it
associated digraph} of $A$, denoted by $D(A)$, has vertex set $V
=\{1$, $2$, $\cdots$, $n\}$ and arc set $E=\{(i, \ j)|\, a_{ij}\neq 0
\}$.
\end{definition}

For a digraph $S$, the {\it associated $(0,1)$-Boolean matrix} of $S$, denoted by $A_{S}$, is a square $(0,1)$-Boolean matrix $A$ with $D(A)$=$S$.
Denoted by $exp(S)$ the primitive $exponent$ of a primitive digraph $S$.
For a primitive $(0,1)$-Boolean $A$, we denote by $exp(A)$ its primitive $exponent$. From combinatorial matrix theory, we know that for a primitive $(0,1)$-Boolean $A$, $exp(A) = exp(D(A))$, and for a primitive digraph $S$, $exp(S)=exp(A_{S})$. As a result, the digraphs can be used to study the primitivity of the $(0,1)$-Boolean matrices.

Primitivity of a square $(0,1)$-Boolean matrix
is of great significance, which is
closely related to many other problems in various areas of pure and
applied mathematics (see \cite{5}, \cite{Kim}, \cite{LB1}, \cite{39},
\cite{55}). In practice, we consider the {\it memoryless
communication system} \cite{39} in communication field, which is
depicted as a digraph $D$ of order $n$. Suppose that each vertex
controls a different piece of information at time $t=0$. At the time
$t=1$, for any vertex $v\in V(D)$, the information $\sigma_{v}$
controlled by vertex $v$ is communicated to each vertex in
$N^{+}(v)$ ({\it out-neighbors}), but $v$ ``forgets" the
information. Proceeding like this, we can ensure that each vertex in
$D$ receives the $n$ pieces of different information at sometime $t$
if $D$ is primitive, where the least time $t$ such that each vertex
in $D$ receive the $n$ pieces of different information is equal to
the primitive exponent of $D$. So, studying the primitvity
of the digraphs is very useful in
information communication field.

It is known that an important topic about the primitivity of the Boolean matrices is to determine the primitivity of a given class of Boolean matrices including the extremal exponents and exponents distribution. Some interesting results about this topic have been shown in many references (see \cite{LB1}, for instance).

Because there is no ambiguity in this paper, for convenience, a directed walk, a directed
path and a directed cycle is abbreviated into a walk, a path and a cycle, respectively. The $girth$ of a strongly connnected digraph, denoted by $g$, is the length of the shortest cycle in this digraph. The $girth$ of a primitive Boolean matrix is defined to be the $girth$ of its associated digraph. In this paper, among all primitive Boolean matrices of order $n$, we consider the primitive exponents of those of girth $g$. For the primitive matrices of both order $n\geq 10$ and girth $g>\frac{n^{2}-4n}{4(n-3)}$, the matrices with primitive exponents in $[2n-2 +(g- 1)(n-3), n+g(n-2)]$ are completely characterized.
\section{Preliminaries}

~~~~In this section, some notations and working lemmas are introduced.

\begin{definition} \label{de2.1}
Let $\{s_{1}$, $s_{2}$, $\cdots$, $s_{\lambda} \}$ be a set of
distinct positive integers with
 gcd$(s_{1}$, $s_{2}$, $\cdots$, $s_{\lambda})$ = 1. The Frobenius number of $
s_{1}$, $s_{2}$, $\cdots$, $s_{\lambda} $, denoted by $\phi(s_{1},
s_{2}, \cdots, s_{\lambda})$, is the smallest nonnegative integer $m$
such that for any nonnegative integers $k\geq m$, there are nonnegative
integers $a_{i} \ (i = 1, 2, \cdots, \lambda)$ such that $k=\sum
\limits_{i = 1}^ \lambda$ $ a_{i}s_{i}$.
\end{definition}

It is well known that
if gcd$(s_{1}, s_{2}) = 1$, then $\phi(s_{1}, s_{2}) = (s_{1} -
1)(s_{2} - 1)$ (see \cite{LB1}, for examlpe).
From Definition \ref {de2.1}, it is easy to see that if there exist
$s_{i}, s_{j}\in\{ s_{1}$, $s_{2}$, $\cdots$, $s_{\lambda} \}$ such
that gcd$(s_{i}, s_{j}) = 1$, then $\phi(s_{1},
s_{2},\cdots,s_{\lambda})\leq \phi(s_{i}, s_{j})$. If $\min \{s_{i}|\, 1\leq i\leq
\lambda\} = 1$, then $\phi(s_{1}, s_{2},\cdots,s_{\lambda})=0$.

In a primitive digraph $S$, the distance from $v_{i}$ to $v_{j}$, denoted by $d(v_{i}, v_{j})$ or $d_{S}(v_{i}, v_{j})$, is the length of the shortest path from $v_{i}$ to $v_{j}$. We denote by $C_{k}$ or {\it $k$-cycle} a cycle with length $k$, and denote by $C(S)$ the cycle length set of digraph $S$. Suppose $C(S)=\{p_{1}$, $p_{2}$,
$\ldots$, $p_{u}\}$. Let $d_{C(S)}(v_{i}$, $v_{j})$ denote the
length of the shortest walk from $v_{i}$ to $v_{j}$ which meets at
least one $p_{i}$-cycle for each $i$, where $i=1$, $2, \cdots, u$. Such a
shortest directed walk is called a $C(S)$-walk from $v_{i}$ to
$v_{j}$. Further, let $d(C(S))=\max\{d_{C(S)}(v_{i}$,
$v_{j})|\, v_{i}$, $v_{j}\in V(S)\}$.

\begin{lemma}{\bf \cite{BLiu}} \label{le2.2}
Let $S$ be a primitive digraph of order $n$
  and $C(S)=\{p_{1}$, $p_{2}$, $\ldots$, $p_{u}\}$.
  Then we have $exp(S)\leq d(C(S))+\phi(p_{1},p_{2},\ldots,p_{u})$.
\end{lemma}

The
union of digraphs $H$ and $G$ is the digraph $G\cup H$ with
vertex set $V(G)\cup V(H)$
 and arc set $E(G)\cup E(H)$. The
$standard$ $n$-$cycle$ is defined to be $\mathscr{C}=(v_{n}$, $v_{n-1}$, $\cdots$ , $v_{2}$,
$v_{1}$, $v_{n})$. Let $D_1=\mathscr{C}\cup(v_{1}$, $v_{n-1})$ and $D_{2}$ =
$D_{1}\cup (v_{2}$, $v_{n})$.

\begin{lemma}{\bf \cite{5}} \label{le2.3}
Let $S$ be a primitive digraph of both order $n$ and girth $g$. Then $exp(S)\leq n+g(n-2)$.
\end{lemma}

\begin{lemma}{\bf \cite{HWL}} \label{le2.4}
Let $S$ be a primitive digraph of order $n$.

$\mathrm{(i)}$ $exp(S)= (n-1)^{2}+1$ if and only if $S\cong D_{1}$;

$\mathrm{(ii)}$ $exp(S)= (n-1)^{2}$ if and only if $S\cong D_{2}$.
\end{lemma}

Lemma \ref{le2.3} gives an upper bound for the primitive exponent of a digraph with both order $n$ and girth $g$. Because a digraph with both order $n$ and girth $n-1$ is isomorphic to $D_{1}$ or $D_{2}$,  from Lemma \ref{le2.4}, we see that the digraphs with both order $n$ and girth $n-1$ are completely characterized.  But for those with girth less than $n-1$, this problem is unsolved as far.

\begin{lemma}{\bf \cite{MiaoZhang}} \label{le2.5}
Let $S$ be a primitive digraph of order $n$ and $|C(S)|\geq 3$. Then
$exp(S)\leq \lfloor\frac{(n-2)^{2}}{2}\rfloor+n$.
\end{lemma}

From Lemma \ref{le2.5}, we get the following corollary.

\begin{corollary}\label{cl2.1}
Let $S$ be a primitive digraph of order $n$. If $exp(S)> \lfloor\frac{(n-2)^{2}}{2}\rfloor+n$, then $|C(S)|= 2$.
\end{corollary}

\begin{lemma}{\bf \cite{SKK}} \label{le2.6}
Let $S$ be a primitive directed graph on $n$ vertices having
cycles of just two lengths, $g$ and $q$, where without loss of generality we take
$g\leq q$. Then $exp(S) \leq 2n -g - 1 + (g- 1)(q - 1)$.
\end{lemma}

\section{Main results}

\begin{lemma}\label{th3.2}
Let $S$ be a primitive digraph of order $n\geq 6$. If $C(S)=\{g, q\}$ where $1\leq g< q\leq n-1$, then $exp(S)\leq 2n-2 +(g- 1)(n-3)$.
\end{lemma}

\begin{proof}
Note that $2n - g- 1 + (g- 1)(q - 1)=2n-2 +(g- 1)(q - 2)$, and $f(q)=2n-2 +(g- 1)(q - 2)$ is  monotone nondecreacing with respect to $q$. By Lemma \ref{le2.6}, the result follows as desired.  $\ \ \ \ \Box$

\end{proof}

Let $n$, $g$ be two positive integers satisfying $\gcd(n, g)=1$, and let $t=\min\{n-g+1, g\}$, $F=\{1, 2, 3, \ldots, t\}$. Let $N\subseteq F$, $D_{g, N}=\mathscr{C}\cup_{i\in N}(v_{i}$,
$v_{g+i-1})$ (see Fig. 3.1).

\vspace{5mm}
\unitlength 1mm \linethickness{0.4pt}
\begin{picture}(94.33,37.33)
\bezier{156}(41.33,20.00)(39.33,6.67)(65.00,7.33)
\bezier{172}(41.33,20.00)(39.33,36.67)(65.00,35.00)
\bezier{176}(65.33,35.00)(93.33,36.33)(91.33,20.33)
\bezier{168}(91.33,20.33)(92.33,5.67)(65.33,7.33)
\put(61.00,7.33){\line(5,3){30.00}}
\put(84.33,9.33){\line(-2,3){17.11}}
\put(49.33,8.67){\circle*{1.49}} \put(91.00,15.33){\circle*{1.49}}
\put(61.33,7.33){\circle*{1.49}} \put(84.33,9.33){\circle*{1.49}}
\put(91.33,25.33){\circle*{1.33}} \put(67.67,35.00){\circle*{1.49}}
\put(41.67,15.33){\circle*{1.49}} \put(41.33,25.67){\circle*{1.33}}
\put(48.00,9.33){\vector(-2,1){3.00}}
\put(41.33,15.67){\vector(0,1){5.67}}
\put(42.33,28.00){\vector(3,4){1.75}}
\put(52.67,34.67){\vector(1,0){5.33}}
\put(68.00,35.00){\vector(1,0){6.33}}
\put(91.33,25.33){\vector(0,-1){6.67}}
\put(86.00,32.00){\vector(4,-3){3.00}}
\put(90.33,14.33){\vector(-2,-3){1.78}}
\put(74.33,7.33){\vector(-1,0){4.33}}
\put(60.67,7.33){\vector(-1,0){4.67}}
\put(49.33,8.67){\line(6,1){41.67}}
\put(49.33,9.00){\vector(1,0){3.00}}
\put(81.33,19.33){\vector(3,2){4.67}}
\put(77.67,19.67){\vector(-2,3){2.00}}
\put(47.67,6.33){\makebox(0,0)[cc]{$v_{1}$}}
\put(39.00,15.00){\makebox(0,0)[cc]{$v_{n}$}}
\put(37.00,25.67){\makebox(0,0)[cc]{$v_{n-1}$}}
\put(61.33,4.00){\makebox(0,0)[cc]{$v_{2}$}}
\put(85.33,7.00){\makebox(0,0)[cc]{$v_{i}$}}
\put(94.50,14.00){\makebox(0,0)[cc]{$v_{g}$}}
\put(96.83,26.00){\makebox(0,0)[cc]{$v_{g+1}$}}
\put(67.33,38.33){\makebox(0,0)[cc]{$v_{g+i-1}$}}
\put(65.67,-2.00){\makebox(0,0)[cc]{Fig. 3.1. $D_{g, N}$}}
\end{picture}

\begin{theorem}\label{th3.3}
Let $r=\max\{u\, |\, u\in N\}$. Then $\mbox{exp}(D_{g,
N})=(n-2)g+1-r+n$.

\end{theorem}

\begin{proof}
{\bf  Case 1}  $r<n-g+1$. Then $g+r-1< n$ and $d(C(D_{g, N}))=d_{C(D_{g,
N})}(v_{n},v_{g+r})=2n-g-r.$ By Lemma \ref{le2.2}, we get
$$\mbox{exp}(D_{g, N})\leq
d(C(D_{g, N}))+\phi(n, g)=2n-g-r+(n-1)(g-1)=(n-2)g+1-r+n.\ $$

Next, we prove $\mbox{exp}(D_{g, N})=2n-g-r+(n-1)(g-1)$. To prove this, we prove that there is no walk of length
$2n-g-r+(n-1)(g-1)-1$ from $v_{n}$ to $v_{g+r}$. Otherwise,
suppose $W$ is a walk of length $2n-g-r+(n-1)(g-1)-1$ from
$v_{n}$ to $v_{g+r}$. Let $\mathcal {P}$ denote the path from $v_{n}$ to
$v_{g+r}$ on cycle $\mathscr{C}$. Denote by $L(\mathcal {P})$ the length of $\mathcal {P}$. Then $L(\mathcal {P})= d(v_{n},
v_{g+r})=n-g-r$ and $\mathcal {P}$ meet only $n$-cycle $\mathscr{C}$ but not any
$g$-cycle. $W$ must contain $\mathcal {P}\cup \mathscr{C}$, some
$g$-cycles and some times of $\mathscr{C}$, that is,
$$2n-g-r+(n-1)(g-1)-1=n-g-r+n+a_{1}n+a_{2}g\ \ (a_{j}\geq
0,\ j=1,2).$$ This induces $\phi(n,g)-1=a_{1}n+a_{2}g\ \ (a_{j}\geq 0,\
j=1,2)$, which contradicts the definition of $\phi(n,g)$. As a result, there is no walk of length
$2n-g-r+(n-1)(g-1)-1$ from $v_{n}$ to $v_{g+r}$. Then
$\mbox{exp}(D_{g, N})=2n-g-r+(n-1)(g-1)$ follows.

{\bf Case 2}  $r=n-g+1$. Then $g+r-1=n$ and $d(C(D_{g, N}))=d_{C(D_{g,
N})}(v_{n},v_{1})=n-1.$ Analogous to Case 1, we get
$\mbox{exp}(D_{g, N})=
(n-1)g=(n-2)g+1-r+n.$ $\ \ \ \Box$
\end{proof}

Let $n$, $g$ ($n\geq 2g$) be two positive integers satisfying $\gcd(n, g)=1$. Let $H$ consist of $n$-cycle $(v_{1}$, $v_{2}$, $\ldots$, $v_{n}$, $v_{1})$ and two $g$-cycles where the two $g$-cycles have no common vertex (see Fig. 3.2). Denote by $C^{1}=(v_{1}$, $v_{2}$, $\ldots$, $v_{g}$, $v_{1})$ and $C^{2}=(v_{k}$, $v_{k+1}$, $\ldots$, $v_{k+g-1}$, $v_{k})$ the two $g$-cycles in $H$ where $k\geq g+1$, $g+k-1\leq n$.

\setlength{\unitlength}{0.7pt}
\begin{center}
\begin{picture}(200,206)
\put(37,104){\circle*{4}}
\put(35,137){\circle*{4}}
\put(167,64){\circle*{4}}
\qbezier(34,120)(34,87)(56,64)\qbezier(56,64)(79,42)(112,42)\qbezier(112,42)(144,42)(167,64)
\qbezier(167,64)(190,87)(190,120)\qbezier(34,120)(34,152)(56,175)\qbezier(56,175)(79,198)(112,198)
\qbezier(112,198)(144,198)(167,175)\qbezier(167,175)(190,152)(190,120)
\put(116,42){\circle*{4}}
\put(160,181){\circle*{4}}
\put(56,65){\circle*{4}}
\qbezier(167,64)(102,84)(37,104)
\put(63,182){\circle*{4}}
\put(188,125){\circle*{4}}
\qbezier(63,182)(125,154)(188,125)
\put(15,105){$v_{1}$}
\put(168,53){$v_{g}$}
\put(13,142){$v_{n}$}
\put(37,58){$v_{2}$}
\put(193,120){$v_{k}$}
\put(25,192){$v_{k+g-1}$}
\qbezier(39,145)(41,159)(44,172)
\put(70,-9){Fig. 3.2. $H$}
\qbezier(39,145)(47,156)(55,166)
\qbezier(188,106)(180,95)(173,84)
\qbezier(188,106)(188,94)(188,81)
\qbezier(119,197)(129,192)(139,187)
\qbezier(119,197)(131,199)(143,200)
\qbezier(140,146)(128,148)(117,150)
\qbezier(140,146)(132,155)(125,163)
\qbezier(68,63)(74,57)(81,50)
\qbezier(81,50)(71,52)(62,54)
\qbezier(110,82)(118,76)(127,70)
\qbezier(110,82)(121,82)(132,81)
\end{picture}
\end{center}

\begin{lemma}\label{le3.4}
$\mbox{exp}(H)\leq(n-1)g+n-2g$.

\end{lemma}

\begin{proof}
To prove this lemma, we consider the following two cases.

{\bf  Case 1}  $d(v_{g+k-1},v_{1})= 1$ and $d(v_{g},v_{k})= 1$. Now, $n=2g$, $k=g+1$, $g+k-1=n$, and $d(C(H))=d_{C(H)}(v_{1},v_{n})=n-1.$ By Lemma \ref{le2.2}, we get
$$\mbox{exp}(H)\leq
d(C(H))+\phi(n, g)=n-1+(n-1)(g-1)=(n-1)g.\ $$

{\bf  Case 2}  $\max\{d(v_{g+k-1},v_{1})$, $d(v_{g},v_{k})\}\geq 2$. Then $2g+1\leq n$. If $d(v_{g+k-1},v_{1})=\max\{d(v_{g+k-1},v_{1})$, $d(v_{g},v_{k})\}$. Then $g+k\leq n$, and $d(C(H))=d_{C(H)}(v_{g+k},v_{n})=n-g-k+n=2n-g-k.$ By Lemma \ref{le2.2}, we get
$$\mbox{exp}(H)\leq
(n-1)(g-1)+2n-g-k=(n-1)g+n-g-k+1.$$
If $d(v_{g},v_{k})=\max\{d(v_{g+k-1},v_{1})$, $d(v_{g},v_{k})\}$, then $k-g\geq2$ and $d(C(H))=d_{C(H)}(v_{g+1},v_{k-1})=k-g-2+n.$ In a same way, we get
$$\mbox{exp}(H)\leq
(n-1)(g-1)+k-g-2+n=(n-1)g+k-g-1.$$ Note that $k-g-1\leq n-2g$, $n-k-g+1\leq n-2g$.
Then the result follows.
$\ \ \ \Box$
\end{proof}

Let $D^{r}=\{D_{g,
N}\, |\, N\subseteq F, \max\{u\, |\, u\in N\}=r\}$. From Theorem \ref{th3.3}, we see that for every digraph $S\in D^{r}$, we have $\mbox{exp}(S)=(n-2)g+1-r+n$.

\begin{theorem}\label{th3.6}
 Suppose $\mathrm{gcd}(g, n) = 1$. Let $S$ be a primitive digraph of both order $n\geq 10$ and girth $g>\frac{n^{2}-4n}{4(n-3)}$, and let $2n-2 +(g- 1)(n-3)< w\leq n+g(n-2)$, $z=(n-2)g+1+n-w$. If $exp(S)=w$, then $S$ is isomorphic to one in $D^{z}$.

\end{theorem}

\begin{proof}
Note that $2n-1 +(g- 1)(n-3)>\lfloor\frac{(n-2)^{2}}{2}\rfloor+n$. By Corollary \ref{cl2.1}, it follows that $C(S)=\{g, q\}$. Note that by Lemma \ref{th3.2}, if $q\leq n-1$, then $\mbox{exp}(S)\leq 2n-2 +(g- 1)(n-3)$. Hence, it follows that $q=n$. It is known that for a primitive diagraph $G$ and its a primitive spanning subgraph $D$, it follows that $exp(G)\leq exp(D)$. Note that $\frac{n^{2}-4n}{4(n-3)}\geq 2$. So, $g\geq 3$. Note that if $n\geq 2g$, then $(n-1)g+n-2g< 2n-1 +(g- 1)(n-3)$. This means that $S$ has no subgraph isomorphic to $H$. Then the result follows from Theorem \ref{th3.3}. $\ \ \ \Box$
\end{proof}

Noting the relation between Boolean matrices and digraphs, we have the following corollary.

\begin{corollary}\label{th3.7}
 Suppose $\mathrm{gcd}(g, n) = 1$. Let $A$ be a primitive Boolean matrix of both order $n\geq 10$ and girth $g>\frac{n^{2}-4n}{4(n-3)}$, and let $2n-2 +(g- 1)(n-3)< w\leq n+g(n-2)$, $z=n+1+g(n-2)-w$. If $exp(A)=w$, then $D(A)$ is isomorphic to one in $D^{z}$.

\end{corollary}

Let $Q_{1}= D_{g, N}$ for $N=\{1\}$, $Q_{2}= D_{g, N}$ for $N=\{1, 2\}$. Note that $D^{1}= \{Q_{1}\}$ and $D^{2}= \{Q_{2}\}$. Then we have the following corollary.

\begin{corollary}\label{th3.8}
 Suppose $\mathrm{gcd}(g, n) = 1$. Let $A$ be a primitive digraph of both order $n\geq 10$ and girth $g>\frac{n^{2}-4n}{4(n-3)}$. Then

 $\mathrm{(i)}$ $\mbox{exp}(A)=(n-2)g+n$ if and only if $D(A)\cong Q_{1}$;

  $\mathrm{(ii)}$ $\mbox{exp}(A)=(n-2)g+n-1$ if and only if $D(A)\cong Q_{2}$.

\end{corollary}

\noindent{\bf Acknowledgment}

We offer many thanks to the referees for their kind reviews and
helpful suggestions.

\small {

}

\end{document}